\documentclass[12pt]{article}
\usepackage{amsfonts,amssymb,latexsym,enumerate}

\newtheorem{definition}{Definition}[section]
\newtheorem{thm}[definition]{Theorem}

\newtheorem{remark}[definition]{Remark}
\newtheorem{lemma}[definition]{Lemma}
\newtheorem{cor}[definition]{Corollary}
\newtheorem{example}[definition]{Example}

\newcommand{\pfbegin}{\noindent {\bf Proof:\ }}
\newcommand{\pfend}{\hfill {\bf $\Box$}}
\newcommand{\exbegin}{\begin{example} \rm }
\newcommand{\exend}{\hfill{\bf $\Box$} \end{example} }
\newcommand{\rmkbegin}{\begin{remark} \rm }
\newcommand{\rmkend}{\hfill{\bf $\Box$} \end{remark} }

\newcommand{\A}{{\mathbb{A}}}

\newcommand{\N}{{\mathbb{N}}}
\newcommand{\C}{{\mathbb{C}}}

\newcommand{\scs}{\scriptstyle}

\newcommand{\ini}{{\rm in}}

\newcounter{secnum}
\newcommand{\sect}[1]
 {\protect\section{#1}
  \protect\setcounter{secnum}{\value{section}}
  \protect\setcounter{equation}{0}
 \protect\renewcommand{\theequation}{\mbox{\arabic{secnum}.\arabic{equation}}}}
\begin{document}

\title{Algorithms for Determining Birationality of Parametrization of
  Affine Curves}
\author{Hyungju Park \\
  Department of Mathematics and Statistics\\
  Oakland University\\
  Rochester, MI 48309\\
  email: park@oakland.edu} \date{} \maketitle

\begin{abstract}
  Let $k$ be an arbitrary field, and $C$ be a curve in $\A^n$
  defined parametrically by $x_1=f_1(t),\ldots,x_n=f_n(t)$, where
  $f_1,\ldots,f_n\in k[t]$.  A necessary and sufficient condition for
  the two function fields $k(t)$ and $k(f_1,\ldots,f_n)$ to be same is
  developed in terms of zero-dimensionality of a derived ideal in the
  bivariate polynomial ring $k[s,t]$. Since zero-dimensionality of
  such an ideal can be readily determined by a Gr\"{o}bner basis
  computation, this gives an algorithm that determines if the
  parametrization $\psi=(f_1,\ldots,f_n): \A \rightarrow C$ is a
  birational equivalence.  We also develop an algorithm that
  determines if $k[t]$ and $k[f_1,\ldots,f_n]$ are same, by which we
  get an algorithm that determines if the parametrization
  $\psi=(f_1,\ldots,f_n): \A \rightarrow C$ is an isomorphism.  We
  include some computational examples showing the application of these
  algorithms.
\end{abstract}


\sect{Introduction}
For an arbitrary field $k$, let $f_1,\ldots, f_n\in k[t]$.  Consider
the curve $C\subset \A_k^n$ given parametrically by
$$x_1=f_1(t),\ldots,x_n=f_n(t).$$ 
This paper is concerned with the following
two questions on this parametrization of $C$:

\begin{enumerate}
\item[\bf Q1.] Can we effectively determine if the parametrization 
  $$\psi =(f_1,\ldots,f_n): \A \rightarrow C$$
  is a birational equivalence,
  or equivalently $k(f_1,\ldots, f_n)=k(t)$?
\item[\bf Q2.] If the answer to Q1 is yes, can we determine if the
  parametrization $\psi$ is an isomorphism, or equivalently
  $k[f_1,\ldots, f_n]=k[t]$?
\end{enumerate}

For the case of two polynomials $f_1, f_2 \in k[t]$, there are several
results available.  S. Abhyankar and T. Moh~\cite{abhyankar:moh} have
given a necessary condition for $k[f_1,f_2]=k[t]$ in terms of degrees
of $f_1$ and $f_2$:

\begin{thm} {\rm (S. Abhyankar and T. Moh, \cite{abhyankar:moh})}
Let $k$ be an arbitrary field of characteristic $p$ ($p=0$
or $p > 0$). Suppose that $f_1$ and $f_2$ are in $k[t]$ with 
$m:=\deg(f_1)\leq n:=\deg(f_2)$, and that $p$ does not divide
$\gcd(m,n)$. If $k[f_1,f_2]=k[t]$, then $m$ divides $n$.
\end{thm}

A. van den Essen and J. Yu~\cite{essen:yu} introduced the notion of
$D$-resultant $D(s)\in k[s]$ of $f$ and $g$, and gave necessary and
sufficient conditions for $k(f,g)=k(t)$ and $k[f,g]=k[t]$ in terms of
$D(s)$ (Theorem 2.1, \cite{essen:yu}).

\exbegin\label{ex:first}
Consider $f_1(t)=t^3, f_2(t)=t^2+t\in k[t]$. The
function field $k(t)$ is an algebraic extension of its subfield
$k(t^3,t^2+t)$.  To determine if $k(t^3,t^2+t)=k(t)$, denote the
field $k(t^3,t^2+t)$ by $K$.  Then
\begin{eqnarray*}
t^2+t\in K &\Longrightarrow &
(t^2+t)^2=t^4+2t^3+t^2=t^4-t + (t^2+t) + 2t^3\in K \\
&\Longrightarrow &
t^4 -t =t(t^3-1) \in K \quad (\mbox{Since}\ t^2+t, 2t^3\in K)\\
&\Longrightarrow &
t\in K \quad (\mbox{Since}\ t^3-1\in K).
\end{eqnarray*}
This shows $k(t^3,t^2+t)=k(t)$. However, since the Abhyankar-Moh
necessary condition is not satisfied, $k[f_1,f_2]\subsetneq k[t]$.
This means that the parametrization $\psi=(t^3,t^2+t):\A\rightarrow C$
is a birational equivalence, but not an isomorphism.
\exend

In this paper, we give algorithmic solutions to Q1 and Q2 in the
general case of $n$ polynomials using the method of Gr\"{o}bner bases.

\sect{Birational and Isomorphic Parametrization of Curves}
Throughout this paper, unless there is a possibility of
confusion, we will use the shorthand notation $\A^n$ for the
affine space $\A_{k}^n$.

For an arbitrary field $k$, consider the curve $C\subset \A^n$ given
parametrically by 
$$x_1=f_1(t),\, \ldots,\, x_n=f_n(t),$$
where $f_1,\ldots,f_n\in k[t]$.  The
morphism $\psi$ defined by
$$\psi=(f_1,\ldots, f_n):\A \rightarrow \A^n,$$
will be simply
referred to as the parametrization of $C$ by polynomials
$f_1,\ldots,f_n$.  The curve $C$ is the Zariski closure 
$\overline{{\rm Im}(\psi)}$ of ${\rm Im}(\psi)$ in $\A^n$, and we occasionally
identify the morphism $\psi:\A\rightarrow \A^n$ with the induced
morphism $\psi:\A \rightarrow C$,

Problem Q1 is equivalent to determining if the induced map of the
functions fields
\begin{eqnarray*}
\psi^*: \ K(\A^n)=k(x_1,\ldots,x_n) &\longrightarrow& K(\A)=k(t) \\
        x_i &\longmapsto & f_i.
\end{eqnarray*}
is surjective, i.e.
$$k(f_1,\ldots,f_n) = k(t).$$
So we need to understand what conditions
on the polynomials $f_1,\ldots,f_n$ will guarantee $k(f_1,\ldots,f_n)
= k(t)$.

For Problem Q2, consider the induced $k$-algebra homomorphism
$\psi^*$ of the coordinate rings
\begin{eqnarray*}
\psi^*: \ A(\A^n)=k[x_1,\ldots,x_n] &\longrightarrow& A(\A)=k[t] \\
        x_i &\longmapsto & f_i.
\end{eqnarray*}
Note that the coordinate ring of $C=\overline{{\rm Im}(\psi)}$ is
$$A(C)=k[x_1,\ldots,x_n]/{\rm Ker}(\psi^*)
={\rm Im}(\psi^*) \cong k[f_1,\ldots,f_n].$$
Proving that $\psi$ is an immersion is equivalent to proving
that their coordinate rings are isomorphic:
$$k[t]=k[f_1,\ldots,f_n].$$

Now we will describe the injectivity and birationality of 
$$\psi: \A\rightarrow C\subset \A^n$$
in terms of a set of bivariate polynomials derived from $f_i$'s.

\begin{lemma}\label{lemma:finite}
  Suppose that $f_1',\ldots,f_n'\in k[t]$ are not identically zero.
  Then, the morphism
  $$\psi:=(f_1,\ldots,f_n):\A \rightarrow \A^n$$
  is finite.
\end{lemma}

\pfbegin 
We may assume, without loss of generality, that $f_1$ is not
a constant.  Consider
\begin{eqnarray*}
\psi^*: k[x_1,\ldots,x_n] &\longrightarrow& k[t] \\
        x_i &\longmapsto & f_i.
\end{eqnarray*}

Dividing $f_1$ by its leading coefficient if necessary, we may assume
$f_1(t)$ is monic. Then the monic polynomial
$$G(T):=f_1(T)-f_1\in k[f_1,\ldots,f_n](T)$$
gives an integral dependence of $t\in k[t]$ on $k[f_1,\ldots,f_n]$.  
\pfend

\bigskip

For $f_1,\ldots,f_n\in k[t]$, we introduce a new variable $s$ and
consider
$$f_1(t),\ldots, f_n(t), f_1(s),\ldots,f_n(s)\in k[s,t].$$
For each $i=1,\ldots,n$, $t-s$ divides $f_i(t)-f_i(s)$ and there
exists $g_i(s,t) \in k[s,t]$ such that $f_i(t)-f_i(s)=(t-s)g_i(s,t)$.
We will identify the fraction $\frac{f_i(t)-f_i(s)}{t-s}$ with the
polynomial $g_i(s,t) \in k[s,t]$.  It is easy to prove that
$g_i(s,s)=f_i'(s)$.  The following theorem characterizes the
algebraic set $V(g_1,\ldots,g_n)\subset \A^2$.

\begin{thm}\label{thm:disjoint}
  For $f_1,\ldots,f_n\in k[t]$, let 
  $$g_i(s,t) :=\frac{f_i(t)-f_i(s)}{t-s}\in k[s,t],\ \ i=1,\ldots, n.$$
  Then, for the morphism $\psi:=(f_1,\ldots,f_n):\A \rightarrow \A^n$,
\begin{eqnarray}
  V(g_1,\ldots,g_n) &=& A_{\psi}\amalg B_{\psi},
\end{eqnarray}
where
\begin{eqnarray*}
  A_{\psi} &=& \{ (a,b) \mid a\neq b\in k\ \mbox{and}\ 
  \psi(a)=\psi(b)\} \\
  B_{\psi} &=& \{(a,a)\mid a\in V(f_1',\ldots,f_n)\}.
\end{eqnarray*}
\end{thm}

\pfbegin 
Suppose $(a,b)\in A_{\psi}$. Then,
\begin{eqnarray*}
  (a,b)\in A_{\psi} &\Longrightarrow&
  a\neq b\ \mbox{and}\ \psi(a)=\psi(b) \\
  &\Longrightarrow& a\neq b\ \mbox{and}\ f_i(a)=f_i(b),\ 
  \forall i=1,\ldots,n \\
  &\Longrightarrow& g_i(a,b)=0,\ \forall i=1,\ldots,n \\
  &\Longrightarrow& (a,b)\in V(g_1,\ldots,g_n).
\end{eqnarray*}
Hence, $A_{\psi}\subset V(g_1,\ldots,g_n)$.

Suppose $(a,a)\in B_{\psi}$. Then,
\begin{eqnarray*}
  (a,a)\in B_{\psi} &\Longrightarrow&
  f_i'(a)=0,\ \forall i=1,\ldots,n \\
  &\Longrightarrow& g_i(a,a)=0,\ \forall i=1,\ldots,n \\
  &\Longrightarrow& (a,a)\in V(g_1,\ldots,g_n).
\end{eqnarray*}
Hence $B_{\psi}\subset V(g_1,\ldots,g_n)$.  Therefore, $A_{\psi}\amalg
B_{\psi}\subset V(g_1,\ldots,g_n)$.

In order to show $V(g_1,\ldots,g_n) \subset A_{\psi}\amalg B_{\psi}$,
let $(a,b)\in V(g_1,\ldots,g_n)$.

If $a\neq b$, then
\begin{eqnarray*}
&& g_i(a,b)=\frac{f_i(b)-f_i(a)}{b-a}=0,\ \forall i=1,\ldots,n \\
&\Longrightarrow& f_i(b)=f_i(a),\ \forall i=1,\ldots,n \\
&\Longrightarrow& \psi(a)=\psi(b) \\
&\Longrightarrow& (a,b)\in A_{\psi}.
\end{eqnarray*}

If $a=b$, then
\begin{eqnarray*}
&& g_i(a,a)=0,\ \forall i=1,\ldots,n \\
&\Longrightarrow& 
a \in V(g_1(s,s),\ldots,g_n(s,s))= V(f_1'(s),\ldots,f_n'(s))\\
&\Longrightarrow& (a,a)\in B_{\psi}.  
\end{eqnarray*}
\pfend

\rmkbegin \label{rmk:singularity}
The set $A_{\psi}$ in Theorem~\ref{thm:disjoint} describes the
multiple points on the curve $C$, while the set $B_{\psi}$ describes
the ramification points (or branch points) on $C$.  More precisely, if
$(a,b)\in A_{\psi}$, then the point $\psi(a)=\psi(b)$ on $C$ has
multiplicity of at least $2$.  If $(a,a)\in B_{\psi}$, then $\psi(a)$
is a ramification point on $C$. Therefore, if $B_{\psi}=\emptyset$,
then the parametrization $\psi: \A\rightarrow C$ is an \'{e}tale
morphism.
\rmkend

\begin{thm}\label{thm:unimodular} 
  For $f_1,\ldots,f_n\in k[t]$, let
  $$g_i(s,t):=\frac{f_i(t)-f_i(s)}{t-s}\in k[s,t],\ \ i=1,\ldots,n.$$
  Suppose that $f_1',\ldots,f_n'\in k[t]$ are not identically
  zero.  Then, the morphism
  $$\psi:=(f_1,\ldots,f_n):\A \rightarrow \A^n$$
  is a closed immersion if and only if $V(g_1,\ldots,g_n)=\emptyset$.
\end{thm}

\pfbegin 
Suppose that $V(g_1,\ldots,g_n)=\emptyset$ and $C$ is the curve in
$\A^n$ given parametrically by $x_1=f_1(t),\ldots,x_n=f_n(t)$.  By
Lemma~\ref{lemma:finite}, $\psi:\A\rightarrow C$ is finite, and thus
proper. Since
$$A_{\psi}=B_{\psi}=\emptyset,$$
by Remark~\ref{rmk:singularity}, $C$ is nonsingular, and thus normal.
Since each fiber of the finite morphism $\psi:\A\rightarrow C$
contains one point, its degree $[K(\A):K(C)]$ is $1$, that is,
$$k(f_1,\ldots,f_n)=k(t).$$ 
The normality of $C$ implies that 
$$A(C)=k[f_1,\ldots,f_n]\subset k[t]$$
is integrally closed in $K(C)=k(f_1,\ldots,f_n)=k(t)$.
Since the normal ring $k[t]$ is integral over $k[f_1,\ldots,f_n]$, 
we conclude 
$$k[f_1,\ldots,f_n]=k[t].$$
Therefore, the morphism
$\psi:\A\rightarrow \A^n$ induces an isomorphism of $\A$ onto the
closed set $C$ in $\A^n$.

Conversely, if $\psi:\A \rightarrow \A^n$ is a closed immersion,
then $C\cong \A$ is nonsingular. By Remark~\ref{rmk:singularity},
$$A_{\psi}=B_{\psi}=\emptyset.$$
Hence $V(g_1,\ldots,g_n)=\emptyset$.
\pfend

\begin{cor}\label{cor:unimodular} 
  For $f_1,\ldots,f_n\in k[t]$, let
  $$g_i(s,t):=\frac{f_i(t)-f_i(s)}{t-s}\in k[s,t],\ \ i=1,\ldots,n.$$
  Suppose that $f_1',\ldots,f_n'\in k[t]$ are not identically zero.
  Then,
  $$k[f_1,\ldots,f_n]=k[t]$$
  if and only if $V(g_1,\ldots,g_n)=\emptyset$.
\end{cor}

\begin{thm}\label{thm:birational}
  For $f_1,\ldots,f_n\in k[t]$, let
  $$g_i(s,t):=\frac{f_i(t)-f_i(s)}{t-s}\in k[s,t],\ \ i=1,\ldots,n.$$
  Suppose that $f_1',\ldots,f_n'\in k[t]$ are not identically zero and
  $C$ is the curve in $\A^n$ given parametrically by
  $x_1=f_1(t),\ldots,x_n=f_n(t)$.  Then the parametrization
  $$\psi=(f_1,\ldots,f_n)=\A\rightarrow C$$
  is a birational
  equivalence if and only if $V(g_1,\ldots,g_n)$ is a finite set.
\end{thm}

\pfbegin  By Lemma~\ref{lemma:finite},
the morphism $\psi: \A\rightarrow \A^n$ is finite, and
thus proper. Hence ${\rm Im}(\psi)$ is a Zariski closed
set of dimension $1$ in $\A^n$, and is equal to $C$. \\[1ex]
($\Longleftarrow$) 
Let $V(g_1,\ldots,g_n)=\{ (a_1,b_1),\ldots,, (a_l,b_l)\}$.
Define open sets $U\subset \A$ and $V\subset \A^n$ by
$$U= \A - \{ a_1,\ldots,a_l\}, \ \ 
V= C- \{ \psi(a_1),\ldots, \psi(a_l).$$
Then $\psi$ induces an finite injective morphism
$\psi\!\mid_{U}: U\rightarrow V$. Since each fiber of $\psi\!\mid_{U}$
has one point, its degree $[K(U):K(V)]$ is $1$. Therefore,
$\psi=\A\rightarrow C$ is birational.

\medskip

\noindent ($\Longrightarrow$)
Since $\psi:\A \rightarrow C={\rm Im}(\psi)$ is birational, there
exist open sets $U\subset \A$ and $V\subset C$ such that
$\psi$ induces an isomorphism between them.  Since $\A-U$ is a proper
closed subset of $\A$, the irreducibility of $\A$ forces $\dim(\A-U) <
1$, i.e.  $\A-U$ is a finite set. Therefore, the injectivity of $\psi$
fails only at finitely many points of $\A$, and $A_{\psi}$ is a finite
set.  Since at least one of $f_i'$'s is nonzero, $f_1',\ldots,f_n'\in
k[t]$ have at most finitely many zeros.  This means that $B_{\psi}$ is
a finite set. Hence, by Theorem~\ref{thm:disjoint}
$V(g_1,\ldots,g_n)=A_{\psi}\amalg B_{\psi}$ is a finite set.  
\pfend

\begin{cor}\label{cor:birational}
  For $f_1,\ldots,f_n\in k[t]$, let
  $$g_i(s,t):=\frac{f_i(t)-f_i(s)}{t-s}\in k[s,t],\ \ i=1,\ldots,n.$$
  Suppose that $f_1',\ldots,f_n'\in k[t]$ are not identically zero.
  Then,
  $$k(f_1,\ldots,f_n)=k(t)$$
  if and only if $|V(g_1,\ldots,g_n)| < \infty$.
\end{cor}

\exbegin 
Consider $f_1(t)=t^3$ and $f_2(t)=t^2+t\in k[t]$ of
Example~\ref{ex:first}.  Let us compute
$$V(\frac{f_1(t)-f_1(s)}{t-s},\frac{f_2(t)-f_2(s)}{t-s}).$$
We have to solve
\begin{eqnarray*}
\frac{f_1(t)-f_1(s)}{t-s} &=& \frac{t^3-s^3}{t-s}=t^2+ts+s^2=0\\
\frac{f_2(t)-f_2(s)}{t-s} &=& 
\frac{(t^2-s^2)+(t-s)}{t-s}=t+s+1=0.
\end{eqnarray*}
From the second equation, $t=-s-1$. By putting it into the
first equation, we get
$$(-s-1)^2+(-s-1)s+s^2=s^2+s+1=0.$$
If $k=\C$, then
$$
V(\frac{f_1(t)-f_1(s)}{t-s},\frac{f_2(t)-f_2(s)}{t-s}) = \{(\frac{\scs
  -1+i\sqrt{3}}{\scs 2}, \frac{\scs -1-i\sqrt{3}}{\scs 2}),
(\frac{\scs -1-i\sqrt{3}}{\scs 2}, \frac{\scs -1+i\sqrt{3}}{\scs
  2})\}.
$$
It can be easily seen that 
$V(\frac{f_1(t)-f_1(s)}{t-s},\frac{f_2(t)-f_2(s)}{t-s})$ is
finite over an arbitrary field $k$.
Since this set is finite, Theorem~\ref{thm:birational} confirms our
earlier finding $k(t^3,t^2+t)=k(t)$.  But since this set is nonempty,
by Theorem~\ref{thm:unimodular}, we conclude that $\psi=(f,g): \A
\rightarrow \A^2$ is not a closed immersion, i.e.  $k[t^3,t^2+t]
\subsetneq k[t]$ as predicted by the Abhyankar-Moh
result~\cite{abhyankar:moh}.  
\exend

\sect{Algorithms}\label{sect:main}
For $f_1,\ldots,f_n\in k[s,t]$, define $g_1,\ldots, g_n$ by
$g_i(s,t):=\frac{f_i(t)-f_i(s)}{t-s}$, $i=1,\ldots,n$.
The finiteness condition on
$V(g_1,\ldots,g_n)$ in Theorem~\ref{thm:birational} is equivalent to
the zero-dimensionality of the ideal $I:=\langle g_1,\ldots,g_n
\rangle \subset k[s,t]$. Fix a term ordering $\prec$ on the set of
monomials in $k[s,t]$, and let $h_1,\ldots,h_l\in k[s,t]$ be the
reduced Gr\"{o}bner basis of the ideal $I$ w.r.t $\prec$.
For an arbitrary polynomial $f\in k[s,t]$, denote the initial (or
leading) term of $f$ w.r.t. $\prec$ by $\ini(f)$.  Then $I$
is zero dimensional if and only if there exist $i,j\in \{1,\ldots,l\}$
such that $\ini(h_i)=s^p$ and $\ini(h_j)=t^q$ for some $p,q\in \N$
(see \cite[Theorem 2.2.7]{adams:loustaunau} for a proof).  This
produces the following algorithmic solution to Problem Q1:

\bigskip

{\bf Algorithm 1}
\begin{tabbing}
\qquad Input: \qquad \= $f_1,\ldots,f_n\in k[t]$. \\
\qquad Output: \> yes if $k(f_1,\ldots,f_n)=k(t)$, no otherwise.\\[0.8ex]
\qquad {\bf Step 1:} \> For each $i=1,\ldots, n$, compute 
$g_i:=\frac{f_i(t)-f_i(s)}{t-s}\in k[s,t]$. \\
\qquad {\bf Step 2:} \> Compute the reduced Gr\"{o}bner basis 
$G=\{h_1,\ldots,h_l\}$ \\
\> of the ideal $I:=\langle g_1,\ldots,g_n \rangle \subset k[s,t]$.\\
\qquad {\bf Step 3:} \> Output yes if there exist $i,j\in \{1,\ldots,l\}$ 
such that \\ 
\> $\ini(h_i)=s^p$ and $\ini(h_j)=t^q$ for some $p,q\in \N$. \\
\> Output no otherwise.
\end{tabbing}
Suppose that the ground field $k$ is algebraically closed and
$g_1,\ldots,g_n\in k[s,t]$. Then, by Hilbert Nullstellensatz,
$V(g_1,\ldots,g_n)\subset \A^2$ is empty if and only if
$g_1,\ldots,g_n\in k[s,t]$ generate the unit ideal, i.e.  there exist
$h_1,\ldots,h_n\in k[s,t]$ such that
$$h_1g_1+\cdots +h_ng_n=1.$$
This observation together with Theorem~\ref{thm:unimodular} produces
the following algorithmic solution to Problem Q2:

\bigskip

{\bf Algorithm 2} (over an algebraically closed field $k$)

\begin{tabbing}
\qquad Input: \qquad \= $f_1,\ldots,f_n\in k[t]$. \\
\qquad Output: \> yes if $k[f_1,\ldots,f_n]=k[t]$, no otherwise.\\[0.8ex]
\qquad {\bf Step 1:} \> For each $i=1,\ldots, n$, compute 
$g_i:=\frac{f_i(t)-f_i(s)}{t-s}\in k[s,t]$.\\
\qquad {\bf Step 2:} \> Compute the reduced Gr\"{o}bner basis $G$ of 
the ideal \\
\> $I:=\langle g_1,\ldots,g_n \rangle \subset k[s,t]$.\\
\qquad {\bf Step 3:} \> Output yes if $G=\{ 1\}$. Output no otherwise.
\end{tabbing}

\sect{Examples} 
Most of the examples in this section are worked out with the computer
algebra system {\em Singular}~\cite{singular}. 

\exbegin 
Consider the trivial case of the twisted cubic $C\subset \A^3$ given
parametrically by
$$x=t, y=t^2, z=t^3,$$
which is apparently isomorphic to the affine
line $\A$.  Let $f_1=t, f_2=t^2, f_3=t^3$, and $g_i(s,t):=
\frac{f_i(t)-f_i(s)}{t-s}$ for each $i=1,2,3$.  Since 
$$g_1(s,t)= \frac{f_1(t)-f_1(s)}{t-s} = 1,$$
the reduced Gr\"{o}bner
basis of $\{g_1,g_2,,g_3\}$ w.r.t an arbitrary term ordering is $\{
1\}$.  Therefore, Algorithm~2 confirms that the parametrization
$$\psi=(t,t^2,t^3): \A\rightarrow C$$
is an isomorphism.
\exend

\exbegin 
Consider the curve $C\subset \A^2$ given parametrically by
$$x=f_1(t):=2{t}^{8}+{t}^{4}+3t+1, \ 
y=f_2(t):={t}^{4}-2{t}^{2}+2.$$
Then,
\begin{eqnarray*}
  g_1(s,t) &:=& \frac{f_1(t)-f_1(s)}{t-s} \\
  &=& 2({t}^{7} + {t}^{6}s + {t}^{5}{s}^{2} + {t}^{4}{s}^{3} + 
  t^3{s}^{4} + t^2s^5 + ts^6 + {s}^{7}) \\
  & & + ({t}^{3} + {t}^{2}s + t{s}^{2} + {s}^{3}) +3, \\
  g_2(s,t) &:=& \frac{f_2(t)-f_2(s)}{t-s} \\
  &=& ({t}^{3} + {t}^{2}s + t{s}^{2} + {s}^{3}) -2(t+s).
\end{eqnarray*}
Fix the degree reverse lex order $\prec$ on $k[s,t]$ with $s\prec t$.
Then a computation with {\em Singular} shows that the reduced
Gr\"{o}bner basis of $\{g_1,g_2\}$ w.r.t. $\prec$ is
$\{h_1, h_2, h_3, h_4\}$ where
\begin{eqnarray*}
h_1 &=& t^2s+s^3-2s,\\
h_2 &=& t^3+t^2s+ts^2+s^3-2t-2s,\\
h_3 &=& 8ts^4+8s^5-16ts^2-16s^3+18, \\
h_4 &=& 16s^6-48s^4-18t^2-ts+51s^2.
\end{eqnarray*}
Therefore, Algorithm~1 and Algorithm~2 show that the parametrization
$\psi=(f_1,f_2): \A\rightarrow C$ is a birational equivalence, but not
an isomorphism (over an algebraically closed field). One notes that,
although $k[f_1,f_2]\neq k[t]$, the Abhyankar-Moh necessary condition
\cite{abhyankar:moh} for $k[f_1,f_2]=k[t]$ is satisfied since
$\deg(f_1)$ divides $\deg(f_2)$.  Hence this example confirms that
the Abhyankar-Moh condition is a necessary but not a sufficient
condition for $k[f_1,f_2]=k[t]$.
\exend

\exbegin
Consider the curve $C\subset \A^3$ given parametrically by
$$x=f_1(t):=t^{10}+t^4,\ y=f_2(t):= t^8+2t^2, \ z=f_3(t):= t^6-t^4+1.$$
Then, for
$g_i(s,t):=\frac{f_i(t)-f_i(s)}{t-s}$, $i=1,2,3$, the reduced
Gr\"{o}bner basis $G$ of $\{g_1,g_2,g_3\}$ w.r.t. the degree reverse lex
order is 
$$G=\{ t+s \}.$$
Hence, according to Algorithm~1, the parametrization
$$\psi=(f_1,f_2,f_3):\A\rightarrow C$$
 is not a birational
equivalence.  
\exend


\end{document}